\documentclass{amsproc}
\usepackage{amssymb}
\usepackage{graphicx}
\usepackage{amscd}
\usepackage{amsmath}
\theoremstyle{plain}

\newtheorem{lemma}{Lemma}

\newtheorem{proposition}{Proposition}
\newtheorem{remark}{Remark}

\newtheorem{theorem}{Theorem}
\numberwithin{equation}{section}

\newcommand\bess{\begin{eqnarray*}}
\newcommand\eess{\end{eqnarray*}}

\newcommand\beq{\begin{equation}}
\newcommand\eeq{\end{equation}}

\def\Nset{\mbox{I\kern-.21em N}}
\def\RE{{\mbox{\rm I\kern-.21em R}}}
\def\ZZ{{\mbox{\sf Z\kern-.45em Z}}}
\def\vv{\kern.344em{\rule[.18ex]{.075em}{1.32ex}}\kern-.344em}

\def\eps{\varepsilon} \def\f{\varphi} 
   \def\k{\varkappa}
   \def\lm{\lambda} 
   
\def\iy{\infty}  \def\K{\mathcal K}
\def\F{\mathcal F}  \def\H{\mathcal H}  
 
\def\l{\left}  \def\r{\right}
\def\<{\langle} \def\>{\rangle}
\begin{document}
\title[]{The boundary control approach to  inverse spectral theory}
\date{February, 2008}
\author{Sergei Avdonin and Victor Mikhaylov}
\address{Department of Mathematics and Statistics \\
University of Alaska Fairbanks\\
PO Box 756660\\
Fairbanks, AK 99775} \email{saavdonin@alaska.edu,
Victor.Mikhaylov@iecn.u-nancy.fr} 
 \maketitle

\begin{abstract}
We establish connections between different approaches to inverse
spectral problems: the classical Gelfand--Levitan theory, the
Krein method, the Simon theory, the approach proposed by Remling
and the Boundary Control method. We show that the Boundary Control
approach provides simple and physically motivated proofs of the
central results of other theories.  We demonstrate also the
connections between the dynamical and spectral data and derive the
local version of the classical Gelfand--Levitan equations.
\end{abstract}

\section{Introduction}

In this paper we consider the Schr\"odinger operator
\begin{equation}
\label{Schrodinger} H=-\partial _{x}^{2}+q\left( x\right)
\end{equation}
on $L^{2}\left( \mathbb{R}_{+}\right) \,,\mathbb{R}_{+}:=[0,\infty
),$ with a real-valued locally integrable potential $q$ and
Dirichlet boundary condition at $x=0.$ Let $d\rho(\lambda)$ be the
\emph{spectral measure} corresponding to $H,$ and $ m(z)$ be the
(principal or Dirichlet) \emph{Titchmarsh-Weyl}
$m$\emph{-function}.

In this section we give a brief review of five different
approaches to inverse problems for the operator
(\ref{Schrodinger}): the Gelfand--Levitan theory, the Krein
method, the Simon theory, the Remling approach and the Boundary
Control method. In the next section we describe the Boundary
Control method in more detail and establish its connections with
the other approaches.

\subsection{Gelfand--Levitan theory}
Determining the potential $q$ from the spectral  measure is the
main result of the seminal paper by Gelfand and Levitan \cite{GL}.
To formulate the result let us define  the following functions:
\begin{eqnarray}
\sigma(\lambda)=\left\{\begin{array}l
\rho(\lambda)-\frac{2}{3\pi}{\lambda}^{\frac{3}{2}},\quad \lambda\geqslant 0, \label{sigma} \\
\rho(\lambda),\quad \lambda <0,\end{array}\right.\\
F(x,t)=\int_{-\infty}^{\iy}\frac{\sin{\sqrt{\lambda}x}\sin{\sqrt{\lambda}t}}{\lambda}\,d\sigma(\lambda).
\label{fxt}
\end{eqnarray}
Let $\f(x,\lm)$ be a solution to the equation
\begin{eqnarray}
&-\f''+q(x)\f=\lambda \f,\quad x >0,\label{Schro_OP}
\end{eqnarray}
with the Cauchy data
\begin{equation}
\label{Gener_eig} \varphi(0,\lambda)=0,\quad
\varphi'(0,\lambda)=1.
\end{equation}
The so-called transformation operator transforms the solutions of
(\ref{Schro_OP}), (\ref{Gener_eig}) with zero potential to the
functions $\f(x,\lm)$:
\begin{equation}
\label{Dir_trans}
\f(x,\lambda)=\frac{\sin{\sqrt{\lambda}x}}{\sqrt{\lambda}}+\int_0^x
K(x,t)\frac{\sin{\sqrt{\lambda}t}}{\sqrt{\lambda}}\,dt.
\end{equation}
The kernel $K(x,t)$ solves the Goursat problem
\begin{equation}
\label{gursa_dir_tr} \left\{\begin{array}l
K_{tt}(x,t)-K_{xx}(x,t)+q(x)K(x,t)=0,\\
K(x,0)=0,\quad \frac{d}{dx}K(x,x)=\frac{1}{2}q(x).
\end{array}
\right.
\end{equation}
It was proved in \cite{GL} that $K(x,t)$ satisfies also the
integral (Gelfand--Levitan) equation
\begin{equation}
\label{GL_classic}
F(x,t)+K(x,t)+\int_0^xK(x,s)\,F(s,t)\,ds=0\,,\quad 0\leqslant t<x.
\end{equation}
The potential can be recovered form the solution of this equation
by the rule
\begin{equation}
q(x)=2\frac{d}{dx}K(x,x).
\end{equation}

\subsection{The Krein method.}

In the beginning of fifties M. Krein developed an approach (see
\cite{Kr1,Kr2}) to spectral inverse problems for the string
equation which is different from the Gelfand-Levitan theory. Using
the method of directing functionals developed by himself in the
forties, Krein reduced the inverse problem to solving the linear
integral equation. Later this equation was derived by
Blagoveschenskii \cite{Bl71} and independently by Gopinath and
Sondhi \cite{GoSh1,GoSh2}  using the dynamical approach.

\subsection{Simon approach. }
In \cite{BS1} Barry Simon proposed a new approach to inverse
spectral theory which has got a further development in the paper
by Gesztesy and Simon \cite{BSFG2} (see also an excellent survey
paper \cite{FG}). As the data of inverse problem they used the
Titchmarsh--Weyl $m$-function which is known to be in one-to-one
correspondence with the spectral measure. It was shown in
\cite{BS1} that there exists a unique real valued function $A \in
L^{1}_{loc}\left( \mathbb{R}_{+}\right) $ (the\emph{\
}$A-$\emph{amplitude}) such that
\begin{equation}
\label{A_amplitude} m(-k^{2})=-k-\int_{0}^{\infty }A(t )e^{-2t
k}\,dt\,.
\end{equation}%
The absolute convergence of the integral was proved for $q\in
L^{1}\left( \mathbb{R}_{+}\right) $ and $q\in L^{\infty }\left( \mathbb{R}%
_{+}\right) $ in \cite{BSFG2} for sufficiently large $\Re k$.  In
\cite{AMR} this result was extended to a more broad class of
potentials. In general situation one has an asymptotic equality
\begin{equation}
\label{A_func_asymp} m(-k^{2})=-k-\int_{0}^{a}A(t)e^{-2t
k}\,dt+\mathcal{O}(e^{-2ak})
\end{equation}%
(see \cite{BS1,BSFG2} for details).

Simon \cite{BS1} put forward the local approach to solving the
inverse problems (locality means that the $A-$amplitude on $[0,a]$
completely determines $q$ on the same interval and vice versa).
First, based on (\ref{A_func_asymp}) Simon proved the local
version of the Borg--Marchenko uniqueness theorem:
$m_1(-k^{2})-m_2(-k^{2})=\mathcal{O}(e^{-2ak})$ if and only if
$q_1(x)=q_2(x)$ for $x \in [0,a].$

Second, he described how to recover potential by the
$A-$amplitude. If $A(\cdot,x)$ denotes the $A-$amplitude of the
problem on $[x,\iy)$, then this family satisfies the  nonlinear
integro-differential equation
\begin{equation}
\label{A_Simon} \frac{\partial A(t,x)}{\partial x}=\frac{\partial
A(t,x)}{\partial t}+\int_0^t A(s,x)A(t-s,x)\,ds=0
 \,.
\end{equation}
If one solves this equation with the initial condition
$A(t,0)=A(t)$ in the domain $\l\{(x,t): 0 \leq x\leq a,\; 0 \leq t
\leq a-x \r\}$, then the potential on $[0,a]$ is determined by
\begin{equation}
 \lim_{t\downarrow 0} A(t,x)=q(x)\,, \ 0 \leq x\leq a.
\end{equation}

The $A-$amplitude has the explicit representation through the
spectral measure by the formula derived in \cite{BSFG2}:
 \beq \label{asp} A(t)=-2 \lim_{\eps \to 0}
\int_{\mathbb{R}}e^{-\eps \lm}\,\frac{\sin (2t {\sqrt
\lm})}{{\sqrt \lm}}\,d\rho(\lm) \ \ a.e. \eeq Without the Abelian
regularization the integral need not be convergent (even
conditionally) \cite{BSFG2}.

\subsection{Remling approach. }

Motivated by Simon, Remling \cite{R1,R2} proposed another local
approach to inverse spectral problems based on the theory of de
Branges spaces. He introduced the integral operator $\K$ acting in
the space $\F^T:=L^2(0,T),$ \beq \label{rem} (\K f)(x)=\int_0^T
k(x,t)\,f(t)\,dt\,, \eeq where \beq \label{remk}
k(x,t)=\frac{1}{2}[\phi(x-t)-\phi(x+t)]\,, \ \
\phi(x)=\int_0^{|x|/2}A(t)\,dt\,. \eeq Remling proved that given a
function $A \in L^1(0,T)$, there exists a unique $q \in L^1(0,T)$
such that $A$ is the $A-$amplitude of this $q$ if and only if the
operator $I+\K$ is positive definite in  $\F^T.$ The same
positivity condition was proved in \cite{R2} to be necessary and
sufficient for solvability of the equation (\ref{A_Simon}).

He proved the following representation of the $A-$amplitude
through the regularized spectral measure $d\sigma$: \beq
\label{asr} A(t)=-2  \int_{\mathbb{R}}\frac{\sin (2t {\sqrt
\lm})}{{\sqrt \lm}}\,d\sigma(\lm) \eeq with the convergence in the
sense of distributions. This result is strengthened in Theorem 2,
Sec. 2.4, below.

Remling derived two linear integral equations, \beq \label{rem1}
y(x,t)+\int_0^x k(t,s)y(x,s)\,ds=t\,, \eeq \beq \label{rem2}
z(x,t)+\int_0^x k(t,s)z(x,s)\,ds=\psi(t) \,, \eeq where $0 \leq t
\leq x \leq T$ and $\psi(t)=-1-\int_0^t \phi(s)\,ds.$ The
potential $q(x)$ on $[0,T]$ is uniquely determined by any of the
functions $y$ or $z$:
\begin{equation}
q(x)=\frac{\frac{d^2}{dx^2}y(x,x)}{y(x,x)},\quad
q(x)=\frac{\frac{d^2}{dx^2}z(x,x)}{z(x,x)}.
\end{equation}
A remarkable fact is that equations (\ref{rem1}), (\ref{rem2}) are
almost identical to Krein's equation and equations of the Boundary
Control (see Sec. 1.5 below).

\subsection{The Boundary Control method.}

The Boundary Control (BC) method in inverse problems was developed
about two decades ago by M.~Belishev and his colleagues
\cite{B87,BKu,BKa,BBl,ABI}. As well as methods of  Simon and
Remling, the BC method provides the local approach to inverse
problems developing ideas of A.~Blagoveshchenskii \cite{Bl71} who
was a pioneer of the local approach to the 1d wave equation. It is
worth to notice that the papers by Simon, Gesztesy and Remling
(and also by Krein \cite{Kr1,Kr2}) are based on the spectral
approach, and locality is proved there using sophisticated
analytical tools. In the BC method, locality naturally follows
from the finite speed of the wave propagation.

The main idea of the BC method is to study the dynamic
Dirichlet-to-Neumann map $R:u(0,t)\mapsto u_x(0,t)$ for the wave
equation associated with the operator (\ref{Schrodinger}):
\begin{equation}
u_{tt}-u_{xx}+q(x)u=0, \quad x>0,\ t>0,
\end{equation}
with zero initial conditions and the boundary condition
$u(0,t)=f(t)$. Operator $R$ has the form
\begin{equation}
(Rf)(t)=-f'(t)+\int_0^tr(s)f(t-s)\,ds,
\end{equation}
and function $r(t)$ is considered as inverse data. Let us
introduce the operator acting in $L^2(0,T)$:
\begin{equation}
\label{r-c} ({C}^T  f)(t)=f(t)+\int_0^T
[p(2T-t-s)-p(t-s)]f(s)\,ds\,, \ 0<t<T\,,
\end{equation}
where
\begin{equation}
\label{pt} p(t):=\frac{1}{2}\int_0^{|t|} r\left(s\right)\,ds.
\end{equation}
It is proved (see, e.g. \cite{ABI}) that one can recover the
potential using the unique solution to any of the equations
\begin{eqnarray} C^Tf^T_0(t)=T-t\,, \ t \in [0,T],\label{c1}\\
(C^Tf_1^T)(t)=-((R^T)^*\k^T)(t) \,, \ t \in [0,T]\label{c2},
\end{eqnarray}
where the operator $R^T$ in the second equation is determined by
$r(t)$, $t\in [0,T]$ and $\varkappa(t)=T-t$. Then
\begin{equation}
q(T)=\frac{\frac{d^2}{dT^2}f_j^T(+0)}{f_j^T(+0)},\quad j=0,1.
\end{equation}

It is important to note that the Krein equation, the Remling
equation and the equation of the BC method can be reduced to each
other by simple changes of variables. More exactly,  Krein in
\cite{Kr1,Kr2} considered the problem with Neumann boundary
conditions at $x=0,$ and one of the equations derived in \cite{R1}
can be reduced to the Krein equation. Equations (\ref{c1}),
(\ref{c2}) and (\ref{rem1}), (\ref{rem2}) concerning Dirichlet
conditions can be easily transformed to each other.

The main goal of this paper is to demonstrate the connections
between all approaches mentioned above.  We provide a new proof of
the Gelfand--Levitan equations which demonstrates their local
character. We describe in detail relations between dynamical and
spectral approaches, in particular, we prove convergence a.e. in
formula (\ref{asr}).

\section{The Boundary Control approach. }

The BC method uses the deep connection between inverse problems of
mathematical physics, functional analysis and control theory for
partial differential equations and offers an interesting and
powerful alternative to previous identification techniques based
on spectral or scattering methods. This approach has several
advantages, namely: (i) it maintains linearity (does not introduce
spurious nonlinearities); (ii) it is applicable to a wide range of
linear point and distributed systems and reconstruction
situations; (iii) it can identify coefficients occurring in
highest order terms; (iv) it is, in principle,
dimension-independent; and, finally, (v) it lends itself to
straightforward algorithmic implementations. Being originally
proposed for solving the boundary inverse problem for the
multidimensional wave equation, the BC method has been
successfully applied to all main types of linear equations of
mathematical physics (see the review papers \cite{B97,B07}, monograph
\cite{KKL} and
references therein). In this paper we use this method in 1d
situation applying it to  inverse problems for the operator
(\ref{Schrodinger}) and demonstrate its connections with the
methods described above. We consider here Dirichlet boundary
condition and note that our approach works also for other boundary
conditions (see, e.g. \cite{ALP} for Neumann condition and
\cite{AP} for a non-self-adjoint condition).

\subsection{The initial boundary value problem, Goursat problem.}

Let us consider the initial boundary value problem for the 1d wave
equation:
\begin{equation}
\label{wave_eqn} \left\{
\begin{array}l
u_{tt}(x,t)-u_{xx}(x,t)+q(x)u(x,t)=0, \quad x>0,\ t>0,\\
u(x,0)=u_t(x,0)=0,\ u(0,t)=f(t).
\end{array}
\right.
\end{equation}
Here $q\in L_{loc}^{1}\left( \mathbb{R}_{+}\right)$ and
 $f$ is an arbitrary $L^2_{loc}\left( \mathbb{R}_{+}\right) $
function referred to as a \emph{boundary control}. The solution
$u^f(x,t)$ of the problem (\ref{wave_eqn}) can be written in terms
of the integral kernel $w(x,s)$ which is the unique solution to
the Goursat problem:
\begin{equation}
\label{gursa} \left\{
\begin{array}l
w_{tt}(x,t)-w_{xx}(x,t)+q(x)w(x,t)=0, \quad 0<x<t,\\
w(0,t)=0,\ w(x,x)=-1/2\int_0^xq(s)\,ds.
\end{array}
\right.
\end{equation}
The properties of the solution to the Goursat problem are given in
Appendix.

Let us consider now the dynamical system $(\ref{wave_eqn})$ on the
time interval $[0,T]$ for some $T>0$. 
\begin{proposition}
\label{Prop_wave_sol}
\begin{itemize}
\item[a)] If $q\in C^1(\mathbb{R}_+)$, $f\in C^2(\mathbb{R}_+)$
and $f(0)=f'(0)=0$, then
\begin{equation}
\label{wave_eqn_sol} u^f(x,t)=\left\{\begin{array}l
f(t-x)+\int_x^tw(x,s)f(t-s)\,ds, \quad x \leq t,\\
0, \quad x > t.\end{array}\right .
\end{equation}
is a classical solution to $(\ref{wave_eqn})$. \item[b)] If $q\in
L^1_{loc}(\mathbb{R}_+)$ and $f\in L^2(0,T)$, then formula
$(\ref{wave_eqn_sol})$ represents a unique generalized solution to
the initial-boundary value problem $(\ref{wave_eqn})$  \\ $u^f \in
C([0,T];\H^T),$ where \bess \mathcal{H}=L^2_{loc}(0,\infty) \ \
\mbox{and} \ \ \mathcal{H}^T := \{u \in \mathcal{H}: \, \mbox{
supp } u \subset [0,T]\, \}. \eess
\end{itemize}
\end{proposition}
First statement of the proposition can be checked by direct
calculations. The proof of the second one follows from
Propositions \ref{Prop_Goursat} and \ref{gursa_conv_prop} (see
Appendix).

\subsection{The main operators of the BC method.}

The {\bf response operator} (the dynamical Dirichlet-to-Neumann
map) $R^T$ for the system (\ref{wave_eqn}) is defined in
$\F^T:=L^2(0,T)$ by
\begin{equation}
\label{Response_sa} (R^Tf)(t)=u_{x}^f(0,t), \ t \in (0,T),
\end{equation}
with the domain $\{ f \in C^2([0,T]):\; f(0)=f'(0)=0\}.$ According
to (\ref{wave_eqn_sol}) it has a representation
\begin{equation}
\label{react_rep} (R^Tf)(t)=-f'(t)+\int_0^tr(s)f(t-s)\,ds,
\end{equation}
where $r(t):=w_x(0,t)$ is called the {\bf response function}.

The response operator $R^T$ is completely determined by the
response function on the interval $[0,T]$, and the dynamical
inverse problem can be formulated as follows. Given $r(t), \; t
\in [0,2T],$ find $q(x), \; x \in [0,T].$

Notice that from (\ref{gursa}) one can derive the formula
\begin{equation}
\label{resp_repr}
r(t)=-\frac{1}{2}q\Bigl(\frac{t}{2}\Bigr)-\frac{1}{2}\int_0^t\
q\Bigl(\frac{t-\zeta}{2}\Bigr) v(\zeta,t)\,d\zeta\,,
\end{equation}
where
$$v(\xi,\eta)=w\Bigl(\frac{\eta-\xi}{2},\frac{\eta+\xi}{2}\Bigl)\,.$$

To solve the dynamical inverse problem by the BC method let us
introduce a couple more operators. Proposition \ref{Prop_wave_sol}
implies in particular that the {\bf control operator} ${W}^T$,
$$
{ W}^T:\F^T \mapsto \H^T, \ { W}^T f=u^f( \cdot,T),
$$
is bounded. From (\ref{wave_eqn_sol}) it follows that
\begin{equation}
\label{W_T_repr} (W^Tf)(x)=f(T-x)+ \int_{x}^{T}
w(x,\tau)f(T-\tau)\,d\tau.
\end{equation}
The next statement claims that the  operator $ W^T $ is boundedly
invertible.

\begin{proposition} \label{pr2}  Let $q \in L_{loc}^1(\mathbb{R}_{+})$ and
$T >0$, then for any function $z \in \mathcal{H}^T$, there exists
a unique control $f \in \mathcal{F}^T$ such that
\begin{equation}\label{eq25}
u^f( x,T) = z(x).
\end{equation}
\end{proposition}

\begin{proof}
 According to (\ref{W_T_repr}), condition (\ref{eq25}) is equivalent to
the following integral Volterra equation of the second kind
\begin{equation} \label{volt}
z(x)=f(T-x) + \int_{x}^{T} w(x,\tau)f(T-\tau)\,d\tau \ \ x \in
(0,T)\,.
\end{equation}
The kernel $ w(x,t) $ is  continuous  and therefore equation
(\ref{volt}) is uniquely solvable, which proves the proposition.
\end{proof}

The {\bf connecting operator} $C^T: \F^T \mapsto \F^T,$ plays a
central role in the BC method. It connects the outer space (the
space of controls) of the dynamical system (\ref{wave_eqn}) with
the inner space (the space of waves) being defined by its bilinear
product:
\begin{equation} \label{Ct} \left<C^Tf,g\right>_{\F^T}=
\left<u^f(\cdot,T),u^g(\cdot,T)\right>_{\H^T}
\end{equation}
In other words,
\begin{equation}
\label{CT}  C^T=(W^T)^*W^T\,,
\end{equation}
and Propositions \ref{Prop_wave_sol}, \ref{pr2} imply that this
operator is positive definite, bounded and boundedly invertible on
$\F^T.$

Let $q_n\in C^\infty(\mathbb{R}_+)$, $n=1,2,\ldots$ and $q_n\to q$
in $L^1_{loc}(\mathbb{R}_+)$. We denote by the $r_n(t)$ the
response function corresponding to $q_n$. Formula
$(\ref{resp_repr})$ and Proposition $\ref{gursa_conv_prop}$ (see Appendix)
yields
\begin{equation}
\label{response_conv} r_n \stackrel{L^1_{loc}}\longrightarrow r,
\quad \text{as $n\to\infty$}.
\end{equation}
Let us denote by $\mathcal{L}(\F^T, \H^T)$ the spaces of bounded
operators acting from $\F^T$ to $\H^T$ and use the notation
$\mathcal{L}(\F^T):=\mathcal{L}(\F^T,\F^T)$. Along with $W^T$ and
$C^T$ we consider operators $W_n^T$ and $C_n^T$ corresponding to
smooth potentials $q_n$, $n=1,2,\ldots$.
\begin{lemma}
\label{cor_oper_conv} Let $W_n^T$, $C_n^T$ be as described above,
then:
\begin{eqnarray}
\|W_n^T - W^T\|_{\mathcal{L}(\F^T, \H^T)}\to 0, \quad \text{as $n\to\infty$},\label{W_T_conv}\\
\|C_n^T -C^T\|_{\mathcal{L}(\F^T)}\to 0, \quad \text{as
$n\to\infty$}\label{C_T_conv}.
\end{eqnarray}
\end{lemma}
\begin{proof}
Let us take arbitrary $f\in \mathcal{F}^T$, then from
(\ref{W_T_repr}) we see that
\begin{equation*}
\|(W-W_n)(f)\|^2_{\mathcal{H}^T}\leqslant
\sup_{0<x<s<T}|w(x,s)-w_n(x,s)|T^2\|f\|^2_{\mathcal{H}^T}.
\end{equation*}
Using $(\ref{gursa_conv})$ (see Appendix), we obtain the first statement of the
lemma. The second statement follows from the first one and the
representation of $C^T$ $(\ref{CT})$.
\end{proof}

The remarkable fact is that $C^T$ can be explicitly expressed
through $R^{2T}$ (or through $r(t),\, t \in [0,2T]$).
\begin{proposition} \label{pr3} For $q \in L_{loc}^1(0,\iy)$ and
$T >0$, operator $C^T$ has the form
\begin{equation}
\label{r-c2} ({C}^T  f)(t)=f(t)+\int_0^T c^T(t,s)f(s)\,ds\,, \
0<t<T\,,
\end{equation}
where
\begin{equation}
\label{c_t} c^T(t,s)=[p(2T-t-s)-p(t-s)].
\end{equation}
 and $p(t)$ is defined in (\ref{pt}).
\end{proposition}
\begin{proof}
For smooth potentials formula (\ref{r-c}) is well known (see, e.g
\cite{ABI}), therefore we give here only a sketch of the proof.
One can easily  check that for $q\in C^\infty(\mathbb{R}_+)$ and
any $f,g\in C_0^\infty(0,T)$ the function $U
(s,t):=\left(u^f(\cdot, s), u^g(\cdot, t)\right)_{\mathcal{H}}$
satisfies the equation $$
U_{tt}-U_{ss}=(R^{T}f(s)g(t)-f(s)(R^Tg)(t)\,, \ \ s,t>0\,,$$ with
the boundary and initial conditions $$U(0,t)=0\,, \ \
U(s,0)=U_t(s,0)=0\,.$$ Using the D'Alambert formula gives
representation (\ref{r-c2}). Making use of the results on the
convergence of operators $(\ref{C_T_conv})$ and response functions
$(\ref{response_conv})$, we can claim that representation
\ref{r-c2} is valid also for $q\in L^1_{loc}(\mathbb{R}_+)$.
\end{proof}

\subsection{The Krein type equations.} First we suppose that $q\in C^\infty(\mathbb{R}_+)$
and consider the Cauchy problem
\begin{equation}
\label{Cauchy_pr} -y''+q(x)y=0,\quad x>0\,; \ \ \  y(0)=\alpha\,,
\ y'(0)=\beta\,.
\end{equation}
Let $f^T$ be a solution of the control problem
\begin{equation}
\label{control_f}
(W^Tf^T)(x)=\left\{\begin{array}l y(x),\ 0<x<T,\\
0,\ x>T.\end{array}\right.
\end{equation}
For any $g\in C_0^\infty(0,T)$ the identity
$$u^g(x,T)=\int_0^T\k^T(t)u^g_{tt}(x,t)\,dt, \ \ \k^T(t):=T-t,$$
is valid, and we have
\begin{eqnarray*}
\left<C^Tf^T,g \right>_{\F^T}=\int_0^T y(x)u^g(x,T)\,dx=\int_0^T
y(x)\int_0^T\k^T(t)u^g_{tt}(x,t)\,dt\,dx\\
=\int_0^T\k^T(t)\left[ y(x)u^g_x(x,T)-y_x(x)u^g(x,T) \right]_{x=0}^{x=T}\,dt\\=
\int_0^T\beta \k^T(t)g(t)-\alpha \k^T(t)(R^Tg)(t)\,dt =\left<\beta
\k^T-\alpha (R^T)^*\k^T,g\right>_{\F^T}.
\end{eqnarray*}
Here $(R^T)^*$ is the operator adjoint to $R^T$ in $\F^T$:
\begin{equation}
((R^T)^*f)(t)=f'(t)+\int_t^Tr(s-t)f(s)\,ds.
\end{equation}
We have used the fact that the solution $u^g(x,t)$ is classical
and $u^g(T,T)=u_x^g(T,T)=0$ (see (\ref{wave_eqn_sol})).

Let us denote by $y_j$, $f_j^T$, $j=0,1$, the functions
corresponding to the cases $\alpha=0$, $\beta=1$ and $\alpha=1$,
$\beta=0,$ respectively. Since $g$ is an arbitrary smooth function,
the functions $f_0^T$ and $f_1^T$ satisfy the equations \beq
\label{Main} (C^Tf_0^T)(t)=T-t\,, \ \
(C^Tf_1^T)(t)=-((R^T)^*\k^T)(t)\,, \ \ t \in [0,T]\,. \eeq Using
(\ref{r-c2}) these equations can be rewritten in more detail: \beq
\label{Main0} f^T_0(t)+\int_0^T c^T(t,s)\,f^T_0(s)\,ds\,=T-t\,, \
t \in [0,T]\,, \eeq \beq \label{Main1} f^T_1(t)+\int_0^T
c^T(t,s)\,f^T_1(s)\,ds\,=1-\int_t^T r(s-t)\,(T-s)\,ds \,, \ t \in
[0,T]\,, \eeq Function $c^T$ is defined in (\ref{c_t}), and from
(\ref{Main0}), (\ref{Main1}) it follows that functions $f^T_j$,
$j=0,1,$ possess additional regularity: $f^T_{j}\in H^1(0,T)$.

Taking into account  $(\ref{C_T_conv})$ and
$(\ref{response_conv})$ we can claim that  equations
$(\ref{Main0})$ and (\ref{Main1}) hold for $q\in
L^1_{loc}(\mathbb{R}_+)$ as well.

Using any of functions $f^T_{j}$ one can easily find the potential $q$ in
the following way. From equation (\ref{wave_eqn_sol}) it follows
that $u^f(t-0,t)=f(+0),$ and in particular, $y_j(T)=f_j^T(+0).$
Let us denote $f_j^T(+0)$ by $\mu_j(T)$. Then \beq \label{mainq}
q(T)=\frac{\mu''_j(T)}{\mu_j(T)}\,. \eeq

Equations (\ref{Main})--(\ref{mainq}) were obtained for a matrix
valued $q$ of a class $C^1$ in \cite{ABI}.

In \cite{AMR} we showed that the Titchmarsh--Weyl $m${-function}
(the spectral Dirichlet-to-Neumann map)  and the response operator
(the dynamical Dirichlet-to-Neumann map) are connected by the
Laplace (or Fourier) transform and established the relation
between the $A-$amplitude and the response function: \beq
\label{ar} A(t)=-2r(2t)\,. \eeq Using this relation it is easy to
check that the positivity condition of Remling's operator $I+\K$
is equivalent to the fact that the operator $C^T$ is positive
definite. Equations (\ref{Main0}), (\ref{Main1}) are reduced by
simple changes of variables to equations (\ref{rem1}),
(\ref{rem2}).

The fact that the positivity of $C^T$ represents the necessary and
sufficient condition of the solvability of the inverse problem
was known in the BC community for a long time.
A.~Blagoveshchenskii \cite{Bl71} in 1971 obtained  the necessary
and sufficient conditions of the solvability of the inverse
problem for the 1d wave equation (with smooth density) which are
equivalent to the positivity of $C^T.$ (Certainly these conditions
were in other terms
--- the BC method and the operator $C^T$ were proposed fifteen years
later). Belishev and Ivanov \cite{BI} considered the two velocity
system with smooth matrix-valued potential. In a particular case
when two velocities are equal, their necessary and sufficient
condition is the positivity of $C^T.$ In \cite{AB} necessary and
sufficient conditions for solvability of a nonselfadjoint inverse
problem with a matrix-valued potential were
formulated in terms of $C^T.$

The equivalent necessary and sufficient conditions for the
solvability of the inverse spectral problem for the string
equation (in the form of positivity of certain integral operator)
were obtained by Krein \cite{Kr1}, \cite{Kr2}.

The method proposed in \cite{BI} works also for non smooth
potentials which leads to the following result.

{\it For given $r \in
L^1(0,2T),$ there exists a unique $q \in L^1(0,T)$ such that $r$
is the response function corresponding to the problem
(\ref{wave_eqn}) with this $q$ if and only if the operator $C^T$
constructed by this $r$ according to (\ref{r-c}) is positive
definite.}

 The fact that $r$ and $q$ belong to the same functional
class is confirmed by formula (\ref{resp_repr}).

\subsection{Spectral representation of $r$ and $c^T$.}

The aim of the present section is to obtain the representation for
$c^T(t,s),$ the kernel of the integral part of the operator $C^T,$
and for the response function $r(t)$ in terms of the spectral
measure of operator (\ref{Schrodinger}).

We consider the Schr\"odinger operator with a real valued
potential $q\in L^1_{loc}(\mathbb{R}_+)$ and Dirichlet boundary
condition at $x=0.$ We remind that $\f(x,\lambda)$ is a solution to the equation
\begin{eqnarray}
\label{Schro_OP_1}
&-\f_{xx}+q\f=\lambda \f,\quad x > 0,
\end{eqnarray}
satisfying the initial conditions
\begin{equation}
\label{Gener_eig_1} \f(0,\lambda)=0,\quad \f'(0,\lambda)=1.
\end{equation}
It is known that there exist a spectral measure $d\rho(\lambda)$,
such that for all $f,g\in L^2(\mathbb{R}_+)$:
\begin{eqnarray}
\int_0^\infty f(x)g(x)\,dx=\int_{-\infty}^\infty
(Ff)(\lambda)(Fg)(\lambda)\,d\rho(\lambda), \label{Fourier_int_1}\\
(Ff)(\lambda)=\int_0^\infty
f(x)\f(x,\lambda)\,dx\label{Fourier_int_2}.
\end{eqnarray}
The so-called inverse transformation operator transforms the
solutions of (\ref{Schro_OP_1}), (\ref{Gener_eig_1}) to the
solutions of the same boundary value problem with $q\equiv 0$,
(cf. (\ref{Dir_trans})):
\begin{equation}
\label{Inv_trans}
\frac{\sin{\sqrt{\lambda}x}}{\sqrt{\lambda}}=\f(x,\lambda)+\int_0^x
L(x,t)\f(t,\lambda)\,dt=:({\bf I}_x+ {\bf L}_x)\f,
\end{equation}
where the kernel $L(x,t)$ satisfy the Goursat problem (see, e.g.
\cite{LB,N}):
\begin{equation}
\label{Inv_trans_gursa} \left\{
\begin{array}l
L_{tt}(x,t)-L_{xx}(x,t)-q(t)L(x,t)=0, \quad 0<t<x,\\
L(x,0)=0,\ \frac{d}{dx}L(x,x)=-\frac{1}{2}q(x).
\end{array}
\right.
\end{equation}
Comparing  (\ref{Inv_trans_gursa}) with (\ref{gursa}) we conclude
that $w(x,t)=L(t,x)$, and thus
\begin{equation}
\label{Inv_trans_wave} \f(s,\lambda)+\int_0^s
w(x,s)\f(x,\lambda)\,dx=\frac{\sin{\sqrt{\lambda}s}}{\sqrt{\lambda}}.
\end{equation}
Let us introduce the functions
\begin{equation}
\Phi_n(s,t)=\int_{-\infty}^n\frac{\sin{\sqrt{\lambda}t}\sin{\sqrt{\lambda}s}}{\lambda}\,d\sigma(\lambda),
\end{equation}
where $\sigma(\lambda)$ is defined in (\ref{sigma}). The fact they
are well-defined follows from the proof of the lemma below. The
following result seems to be  classical, although we have not
been able to find it in the literature for the case of Dirichlet
boundary condition. The case of Neumann boundary condition     is
considered in \cite{LB, N} where the convergence of
corresponding analogues of $\Phi_n$ was proven. We provide the
proof here for the sake of completeness.
\begin{lemma}
\label{Phi_conv} As $n\to\infty,$ the sequence of functions
$\Phi_n(s,t)$ converges uniformly on every bounded set in
$\mathbb{R}^2$ to a continuous function $\Phi(s,t)$ differentiable
outside the diagonal.
\end{lemma}
\begin{proof}
We follow the scheme proposed in \cite{LB}, Lemma 2.2.2. In
\cite{L3} it is shown that the sequence of functions
\begin{equation} \label{Levitan_sf}
\Psi_n(t,s)=\int_{-\infty}^n
\f(t,\lambda)\f(s,\lambda)\,d\rho(\lambda)-\int_0^n
\frac{\sin{\sqrt{\lambda}t}\sin{\sqrt{\lambda}s}}{\lambda}\,d\left(\frac{2}{3\pi}\lambda^{\frac{3}{2}}\right),
\end{equation}
converges uniformly on every bounded set to a differentiable
outside the diagonal function, as $n$ tends to infinity. Applying
operators $({\bf I}_s+{\bf L}_s)({\bf I}_t+{\bf L}_t)$ to (\ref{Levitan_sf}) we have:
\begin{eqnarray}
\label{L_eq}
({\bf I}_s+{\bf L}_s)({\bf I}_t+{\bf L}_t) \Psi_n(t,s)=\Phi_n(s,t) \\
- \int_0^n\left(\int_0^tL(t,\tau)\frac{\sin{\sqrt{\lambda}\tau}}{\sqrt{\lambda}}\,d\tau\right)\frac{\sin{\sqrt{\lambda}s}}{\sqrt{\lambda}}\,
d\left(\frac{2}{3\pi}\lambda^{\frac{3}{2}}\right) \notag\\
- \int_0^n\left(\int_0^sL(s,\tau)\frac{\sin{\sqrt{\lambda}\tau}}{\sqrt{\lambda}}\,d\tau\right)\frac{\sin{\sqrt{\lambda}t}}{\sqrt{\lambda}}\,
d\left(\frac{2}{3\pi}\lambda^{\frac{3}{2}}\right) \notag\\
- \int_0^n\left(\int_0^tL(t,\tau)\frac{\sin{\sqrt{\lambda}\tau}}{\sqrt{\lambda}}\,d\tau\right)
\left(\int_0^sL(s,\tau)\frac{\sin{\sqrt{\lambda}\tau}}{\sqrt{\lambda}}\,d\tau\right)\,d\left(\frac{2}{3\pi}\lambda^{\frac{3}{2}}\right).\notag
\end{eqnarray}
The sum of the last three terms in the right hand side of the
above expression converges to
$-L(s,t)-L(t,s)-\int_0^{\min{\{s,t\}}}L(s,\tau)L(t,\tau)\,d\tau$.
This fact and the convergence of the left hand side of
(\ref{L_eq}) imply the statement of the Lemma.
\end{proof}

The following theorem gives an expression for the  kernel of the
integral part of the operator $C^T$ in terms of the spectral
measure.
\begin{theorem}
\label{Th_Kernel} The kernel $c^T(s,t)$ admits the following
representation:
\begin{equation}
\label{C_T_eqiv} c^T(s,t)=\int_{-\infty}^\infty
\frac{\sin{\sqrt{\lambda}(T-t)}\sin{\sqrt{\lambda}(T-s)}}{\lambda}\,d\sigma(\lambda),\quad
s,t\in [0,T],
\end{equation}
where the integral in the right-hand side of (\ref{C_T_eqiv})
converges uniformly on $[0,T]\times [0,T]$.
\end{theorem}

\begin{proof}
Let  $f,g\in \mathcal{F}^T$. Using
(\ref{Fourier_int_1}) and (\ref{Fourier_int_2}), we rewrite
$\left<C^Tf,g\right>_\mathcal{F^T}$ as
\begin{equation}
\label{C_T_1} \left<C^Tf,g\right>_{\mathcal{F}^T}=\int_0^T
u^f(x,T)u^g(x,T)\,dx=\int_{-\infty}^{\infty}
(Fu^f)(\lambda,T)(Fu^g)(\lambda,T)\,d\rho(\lambda).
\end{equation}
Here (see also (\ref{wave_eqn_sol}))
\begin{eqnarray*}
(Fu^f)(\lambda,T)=\int_0^T \f(x,\lambda)u^f(x,T)\,dx=\\
\int_0^T
\f(x,\lambda)\left(f(T-x)+\int_0^Tw(x,s)f(T-s)\,ds\right)\,dx.
\end{eqnarray*}
Changing the order of integration and using the fact that
$w(x,s)=0$ for $s<x$, we arrive at
\begin{equation}
\label{C_T_2} (Fu^f)(\lambda,T)=\int_0^T
f(T-s)\left(\f(s,\lambda)+\int_0^s
w(x,s)\f(x,\lambda)\,dx\right)\,ds.
\end{equation}
Making use of (\ref{C_T_2}) and (\ref{Inv_trans_wave}), we can
rewrite (\ref{C_T_1}) as
\begin{equation*}
\left<C^Tf,g\right>_{\mathcal{F}^T}=\int_{-\infty}^\infty\,\int_0^T\,\int_0^T\,
\,\frac{\sin{\sqrt{\lambda}(T-t)}\sin{\sqrt{\lambda}(T-s)}}{\lambda}f(t)g(s)\,dt\,ds\,d\rho(\lambda).
\end{equation*}
Comparing the last formula with $(\ref{r-c2})$, we see that
\begin{eqnarray}
\label{C_T_3} \int_{-\infty}^\infty\,\int_0^T\,\int_0^T\,
\,\frac{\sin{\sqrt{\lambda}(T-t)}\sin{\sqrt{\lambda}(T-s)}}{\lambda}f(t)g(s)\,dt\,ds\,d\rho(\lambda)=\\
\int_0^T
f(s)g(s)\,ds+\int_0^T\,\int_0^T\,c^T(s,t)g(s)f(t)\,dt\,ds.\notag
\end{eqnarray}
Now we make use of the $\sin$ transform: for all $h, j\in
L^2(\mathbb{R}_+)$
\begin{eqnarray*}
\widehat h(\lambda)=\int_0^\infty
h(x)\frac{\sin{(\sqrt{\lambda}x)}}{\sqrt{\lambda}}\,dx,\quad
h(x)=\int_0^\infty \widehat h(\lambda)\sin{(\sqrt{\lambda}x)}\,d\left(\frac{2}{3\pi}{\lambda}^\frac{3}{2}\right),\\
\int_0^\infty h(x)j(x)\,dx=\int_0^\infty \widehat
h(\lambda)\widehat
j(\lambda)\,d\left(\frac{2}{3\pi}{\lambda}^\frac{3}{2}\right).
\end{eqnarray*}
Let us extend the functions $f$ and $g$ to the whole real axis
 setting $f(t)=g(t)=0$ for $t>T$ and $t<0$ and use the
notation $f_T(s)=f(T-s)$. Then we can rewrite the first term in
the right hand side of (\ref{C_T_3}) as
\begin{eqnarray}
\label{Delta_tod} \int_0^T f(t)g(t)\,dt=\int_0^\infty
f(T-s)g(T-s)\,ds \\
= \int_0^\infty \widehat f_T(\lambda)\widehat
g_T(\lambda)\,d\left(\frac{2}{3\pi}{\lambda}^\frac{3}{2}\right) \notag\\
= \int_0^\infty\,\int_0^T\, \int_0^T\,
\frac{\sin{\sqrt{\lambda}(T-t)}\sin{\sqrt{\lambda}(T-s)}}{\lambda}f(t)g(s)\,dt\,ds\,d\left(\frac{2}{3\pi}{\lambda}^\frac{3}{2}\right).\notag
\end{eqnarray}
Plugging (\ref{Delta_tod}) in (\ref{C_T_3}), we get
\begin{eqnarray}
\label{C_T_5} \int_{-\infty}^\infty\,\int_0^T\,\int_0^T\,
\,\frac{\sin{\sqrt{\lambda}(T-t)}\sin{\sqrt{\lambda}(T-s)}}{\lambda}f(t)g(s)\,dt\,ds\,d\sigma(\lambda) \\
= \int_0^T\,\int_0^T\,c^T(s,t)f(t)g(s)\,dt\,ds.\notag
\end{eqnarray}
In the last formula the function
$$
C(s,t):=\int_{-\infty}^\infty\,\frac{\sin{\sqrt{\lambda}(T-t)}\sin{\sqrt{\lambda}(T-s)}}{\lambda}\,d\sigma(\lambda)
$$
is a distributional kernel, whose action on functions $f$, $g$ is
defined by the right hand side of (\ref{C_T_5}). On the other hand,
comparing $C(s,t)$ with $\Phi(s,t)$, we see that
$C(s,t)=\Phi(T-s,T-t)$ and according to Lemma \ref{Phi_conv},
$C(s,t)$ is a continuous function on $[0,T]\times [0,T]$. Since
(\ref{C_T_5}) holds for arbitrary $f,g\in \mathcal{F}^T$, we
deduce that
\begin{equation}
\label{C_T_4} c^T(s,t)=C(s,t)=\int_{-\infty}^\infty
\frac{\sin{\sqrt{\lambda}(T-t)}\sin{\sqrt{\lambda}(T-s)}}{\lambda}\,d\sigma(\lambda),\quad
t,s\in [0,T].
\end{equation}
\end{proof}
Using the representation for $c^T(t,s)$ obtained in Theorem
$\ref{Th_Kernel}$, we can derive the formula for the response
function: \begin{theorem} \label{Th_Resp} The representation for
the response function $r$
\begin{equation}
\label{Resp_mes_conn} r(t)=\int_{-\infty}^\infty
\frac{\sin{\sqrt{\lambda}t}}{\sqrt{\lambda}}\,d\sigma(\lambda),\,
\end{equation}
holds for almost all $t \in [0,+\infty)$.
\end{theorem}
\begin{proof}
Let us note that
\begin{equation}
\label{Ph_i}
\Phi(s,t)=\int_{-\infty}^\infty
\frac{\sin{\sqrt{\lambda}t}\sin{\sqrt{\lambda}s}}{\lambda}\,d\sigma(\lambda)=c^T(T-t,T-s),\quad
t,s\in [0,T].
\end{equation}
Using (\ref{c_t}), we have
\begin{equation}
\label{Phi_1}
c^T(T-t,T-s)=\frac{1}{2}\int_{|t-s|}^{t+s}r(\tau)\,d\tau,\quad
t,s\in [0,T].
\end{equation}
The integral in (\ref{Ph_i}) can be rewritten as
\begin{eqnarray}
\Phi(s,t)=\frac{1}{2}\int_{-\infty}^\infty
\frac{(\cos{\sqrt{\lambda}(s+t)}-1-(\cos{\sqrt{\lambda}|s-t|}-1)}{\lambda}\,d\sigma(\lambda)=\label{Phi_2}\\
\frac{1}{2}\int_{-\infty}^\infty\int_{|t-s|}^{t+s}\frac{\sin{\sqrt{\lambda}\theta}}{\sqrt{\lambda}}\,d\theta\,d\sigma(\lambda),\quad
t,s\in [0,T].\notag
\end{eqnarray}
Equating the expressions in  (\ref{Phi_1}) and (\ref{Phi_2})
for $t=s$ we get
\begin{equation}
2c^T(T-t,T-t)=\int_0^{2t}r(\tau)\,d\tau=
\int_{-\infty}^\infty\int_0^{2t}\frac{\sin{\sqrt{\lambda}\theta}}{\sqrt{\lambda}}\,d\theta\,d\sigma(\lambda),\quad
t\in [0,T].
\end{equation}
According to (\ref{react_rep}),  $r\in L^1(0,T)$, so we can use
the Lebesgue theorem and differentiate the last equation. We
obtain the following equality almost everywhere on $[0,2T]$:
\begin{equation}
r(t)=\int_{-\infty}^\infty\frac{\sin{\sqrt{\lambda}t}}{\sqrt{\lambda}}\,d\sigma(\lambda).
\end{equation}
Since the parameter $T$  can be chosen arbitrary
large, the last formula proves the statement of the proposition.
\end{proof}

A direct consequence of this theorem is that integral in formula
(\ref{asr}) converges for almost all $t \in [0,+\infty)$.

The finite speed of the wave propagation (equal to one) in the
equation (\ref{wave_eqn}) implies the local nature of the response
function $r(t)$: the values of $r(t)$, $t\in [0,2T]$ are
determined by the potential $q(x)$, $x\in [0,T]$. Therefore,
 if we are interested in the spectral representation of
$c^T(s,t)$ for $s,t\in [0,T]$ and of $r(t)$ for $t\in [0,2T]$ in
 we can
replace  the formulas (\ref{C_T_eqiv}), (\ref{Resp_mes_conn})
the regularized spectral function
$\sigma(\lambda)$ by (for example) any of the following functions:
\begin{equation}
\sigma_{\operatorname{tr}}(\lambda)=\left\{\begin{array}l
\rho_{\operatorname{tr}}(\lambda)-\frac{2}{3\pi}{\lambda}^{\frac{3}{2}},
\quad \lambda\geqslant 0,\\
\rho_{\operatorname{tr}}(\lambda),\quad \lambda
<0,\end{array}\right.,\quad
\sigma_{d}(\lambda)=\left\{\begin{array}l
\rho_{d}(\lambda)-\rho_0(\lambda),\quad \lambda\geqslant 0,\\
\rho_{d}(\lambda),\quad \lambda <0.\end{array}\right.
\end{equation}
Here $\rho_{tr}$ is the spectral function corresponding to the
truncated potential: $q_T(x)=q(x)$ for $0 \leq x \leq T$ and $q_T(x)=\tilde
q(x)$ for $x>T$ with arbitrary locally integrable $\tilde q$; $\rho_d(\lambda)$ is the
spectral function associated to the discrete problem on the
interval $(0,T)$ with the potential $q_d(x)=q(x)$, $x\in [0,T]$
and $\rho_0(\lambda)$ is the spectral function associated to the
discrete problem on $[0,T]$ with zero potential. (Any
self-adjoint boundary condition can be prescribed at $x=T.$)

\subsection{Gelfand-Levitan equations.}

In this section, using the BC approach we derive the local version
of the classical Gelfand-Levitan equations (\ref{GL_classic}). The
proof is based on the fact that the kernel $K$ of the
transformation operator (\ref{Dir_trans}) satisfies the Goursat
problem (\ref{gursa_dir_tr}). We show that the kernel $v$ of the
operator $(W^T)^{-1}$ (inverse to the control operator
$W^T$) satisfies a similar Goursat problem. We observe that the
operator $(W^T)^{-1}: \mathcal{H}^T\mapsto \mathcal{F}^T$ can be
constructed in the following way. Let us consider the initial-boundary
value problem.
\begin{equation}
\label{wave_eqn_rev} \left\{
\begin{array}l
u_{tt}(x,t)-u_{xx}(x,t)+q(x)u(x,t)=0, \quad 0< x,\, 0< t< T,\\
u(x,T)=a(x);\ \ u(x,t)=0,\ x>t,
\end{array}
\right.
\end{equation}
and denote by $u^a(x,t)$ the solution of this problem. Basing on
the uniqueness of the solutions to the initial boundary value
problems (\ref{wave_eqn}) and (\ref{wave_eqn_rev}) one can check
that
\begin{equation*}
\left(\left(W^T\right)^{-1}a\right)\left(t\right)=u^a(0,t),\quad
0<t<T,
\end{equation*}
(see, e.g. \cite{ABe} for more details). When $q\in
C^1_{loc}(\mathbb{R}_+)$ and $a\in C^1[0,T]$, $a(0)=0$, $u^a(x,t)$
is a classical solution and admits the representation
\begin{equation}
\label{wave_eqn_sol_rev} u^a(x,t)=\left\{\begin{array}l
a(x-t+T)+\int_0^{t-x}v(T-x,s,T-t)a(T-s)\,ds, \quad x\leqslant t,\\
0, \quad x> t, \end{array}\right .
\end{equation}
in terms of the solution $v(x,s,t)$ to the following Goursat
problem:
\begin{equation}
\label{gursa_rev} \left\{
\begin{array}l
v_{tt}(x,s,t)-v_{ss}(x,s,t)+q(T-s)v(x,s,t)=0, \quad 0<s<x-t,\\
v(x,s,0)=0,\, \frac{d}{dt}v(x,t,x-t)=\frac{1}{2}q(T-x+t).
\end{array}
\right.
\end{equation}
By the analogy with Proposition
\ref{Prop_wave_sol} one can show that formula
(\ref{wave_eqn_sol_rev}) gives a generalized solution in the case
of non-smooth potential $q$ and boundary condition $a$. From
representation (\ref{wave_eqn_sol_rev}), the formula for
$(W^T)^{-1}$ immediately follows:
\begin{equation}
\label{W_inv}
\left(\left(W^T\right)^{-1}a\right)(t)=a(T-t)+\int_0^t
V(y,t)a(T-y)\,dy,
\end{equation}
Here the kernel $V(s,t)$ satisfies the Goursat problem
\begin{equation}
\label{gursa_rev_1} \left\{
\begin{array}l
V_{tt}(s,t)-V_{ss}(s,t)+q(T-s)V(s,t)=0, \quad 0<s<t,\\
V(s,T)=0,\ \frac{d}{dt}V(t,t)=\frac{1}{2}q(t).
\end{array}
\right.
\end{equation}
Let us introduce the following operators
\begin{eqnarray}
J_T: L_2(0,T)\mapsto L_2(0,T),\quad (J_T a)(y)=a(T-y),\notag\\
K: L_2(0,T)\mapsto L_2(0,T),\quad (Ka)(t)=\int_0^t
V(y,t)a(y)\,dy,\quad t\in (0,T),\notag\\
K^*: L_2(0,T)\mapsto L_2(0,T),\quad (K^*b)(t)=\int_t^T
V(t,y)b(y)\,dy,\quad t\in (0,T).\notag
\end{eqnarray}
Using these definitions, we can rewrite (\ref{W_inv}) as
\begin{equation}
\left(W^T\right)^{-1}a=(I+K)J_Ta.
\end{equation}
Proposition \ref{pr2} and formula $(W^T)^*=J_T^*(I+K^*)$ yield
\begin{proposition}
\label{Rem_K_adj_inv} The operator $I+K^*: L^2(0,T)\mapsto
L^2(0,T)$ is boundedly invertible.
\end{proposition}
For arbitrary $f,g\in \mathcal{F}^T$, by the definition of $C^T$
we have:
\begin{equation}
\label{C_T_R1}
(C^Tf,g)_{\mathcal{F}^T}=(W^Tf,W^Tg)_{\mathcal{H}^T}.
\end{equation}
Let us put $f=(W^T)^{-1}a$, $g=(W^T)^{-1}b$, $a,b\in
\mathcal{H}^T$ and rewrite (\ref{C_T_R1}) as
\begin{equation}
\label{C_T_R2}
(C^T(I+K)J_Ta,(I+K)J_Tb)_{\mathcal{F}^T}=(a,b)_{\mathcal{H}^T}=(J_Ta,J_Tb)_{\mathcal{H}^T},
\end{equation}
Since (\ref{C_T_R2}) holds for all $a,b\in \mathcal{H}^T$, this
leads to the following operator equation
\begin{equation}
\label{Oper_eqn_inv} (I+K)^*C^T(I+K)=I.
\end{equation}
Introducing the operator
$$
(C_Tf)(t)=\int_0^Tc^T(s,t)f(s)\,ds,
$$
and using (\ref{r-c2}) we can rewrite (\ref{Oper_eqn_inv}) as
\begin{equation}
\label{Oper_eqn_1} K^*+(I+K^*)(K+C_T+C_TK)=0.
\end{equation}

The function $V(y,t)$ was defined in (\ref{gursa_rev_1}) for
$0\leqslant y\leqslant t\leqslant T$, let us continue it by zero
in the domain $t<y\leqslant T$ and introduce the function
$\phi_y(t)$, $y,t\in [0,T]$ by the rule
\begin{equation}
\phi_y(t)=V(y,t)+c^T(y,t)+\int_0^Tc^T(t,s)V(y,s)\,ds.
\end{equation}
The equality (\ref{Oper_eqn_1}) implies
\begin{equation}
V(t,y)+\phi_y(t)+\int_0^T V(t,z)\phi_y(z)\,dz=0,\quad x,t\in
(0,T).
\end{equation}
Since $V(t,y)=0$ for $0<y<t<T$, we obtain that
\begin{equation}
\phi_y(t)+\int_t^T V(t,z)\phi_y(z)\,dz=0,\quad 0<y<t<T.
\end{equation}
Rewriting this equation as
\begin{equation}
((I+K^*)\phi_y)(t)=0,\quad 0<y<t<T,
\end{equation}
and taking into account the invertibility of $I+K^*$ (see
Proposition \ref{Rem_K_adj_inv}), we get
\begin{equation}
\label{G_L}
\phi_y(t)=V(y,t)+c^T(y,t)+\int_y^Tc^T(t,s)V(y,s)\,ds=0,\quad
0<y<t<T.
\end{equation}
Let us formulate this result as
\begin{theorem}
\label{Th_GL} The kernels of operators $C^T$ and  $K$ satisfy the
following integral equation
\begin{equation}
\label{G_LN} V(y,t)+c^T(y,t)+\int_y^Tc^T(t,s)V(y,s)\,ds=0,\quad
0<y<t<T.
\end{equation}
\end{theorem}
Solving the equation (\ref{G_LN}) for all $y \in (0,T)$ we can
recover the potential using
\begin{equation*}
q(y)=2\frac{d}{dx}V(y,y).
\end{equation*}

It is easy to see that the kernel $V$ is connected with the kernel
of the transformation operator (\ref{Dir_trans}) by the rule
$V(T-y,T-t)=K(y,s)$ and $c^T$ is similarly related to $F$ defined
in  (\ref{fxt}): $c^T(T-x,T-t)=F(x,t)$. Therefore, equations
(\ref{G_LN}) can be rewritten in a classical form
(\ref{GL_classic}). On the other hand, equations (\ref{G_LN}) have
clearly a local character since $V(y,t)$ and $c^T(y,t)$ are
completely determined by $q(y)$ on the interval $[0,T].$

\section{Appendix}

The Goursat problem was studied in \cite[Sec. II.4]{TS} for smooth
$q$, but the method works for $q \in L^1(0,a)$ as well (see
\cite{ALP,ALP1,AMR,AK}).
\begin{proposition}
\label{Prop_Goursat}
\begin{itemize}
\item[a)] If $q\in L^1_{loc}(\mathbb{R}_+)$, then the generalized
solution $w(x,s)$ to the Goursat problem $(\ref{gursa})$ is a
continuous function and
\begin{eqnarray}
|w(x,s)|\leqslant
\Bigl(\frac{1}{2}\int_0^{\frac{s+x}{2}}|q(\alpha)|\,d\alpha\Bigr)
\operatorname{exp}\Bigl\{\frac{s-x}{4}\int_0^{\frac{s+x}{2}}|q(\alpha)|\,d\alpha\Bigr\}\label{gursa_est},\\
w_x(\cdot,s), w_s(\cdot,s), w_x(x,\cdot), w_s(x,\cdot)\in
L_{1,\,loc}(\mathbb{R}_+).\label{gursa_deriv}
\end{eqnarray}
Partial derivatives in $(\ref{gursa_deriv})$ continuously in
$L^1_{loc}(\mathbb{R}_+)$ depend on parameters $x$, $s$. The
equation in $(\ref{gursa})$ holds almost everywhere and the
boundary conditions are satisfied in the classical sense.

\item[b)] If $q\in C_{loc}(\mathbb{R}_+)$, then the generalized
solution $ w(x,s)$ to the Goursat problem $(\ref{gursa})$ is
$C^1$-smooth, equation and boundary conditions are satisfied in
the classical sense.

\item[c)] If $q\in C^1_{loc}(\mathbb{R}_+)$, then the solution to
the Goursat problem $(\ref{gursa})$ is classical, all its
derivatives up to the second order are continuous.
\end{itemize}
\end{proposition}
\begin{proof}
By setting $\xi=s-x$, $\eta=s+x$, and
\begin{equation}
\label{Var_change}
v(\xi,\eta)=w\Bigl(\frac{\eta-\xi}{2},\frac{\eta+\xi}{2}\Bigl),
\end{equation}
equation $(\ref{gursa})$ reduces to
\begin{equation}
\label{gursa_trans} \left\{
\begin{array}l
v_{\xi\eta}-\frac{1}{4}q(\frac{\eta-\xi}{2})v=0, \quad 0<\xi<\eta,\\
v(\eta,\eta)=0,\
v(0,\eta)=-\frac{1}{2}\int_0^{\eta/2}q(\alpha)\,d\alpha.
\end{array}
\right.
\end{equation}
The boundary value problem $(\ref{gursa_trans})$ is equivalent to the
integral equation
\begin{equation}
\label{wave_gursa_new_2}
v(\xi,\eta)=-\frac{1}{2}\int_{\xi/2}^{\eta/2}q(\alpha)\,d\alpha-\frac{1}{4}
\int_0^\xi\,d\xi_1\int_\xi^\eta\,d\eta_1q\Bigl(\frac{\eta_1-\xi_1}{2}\Bigr)v(\xi_1,\eta_1).
\end{equation}
We introduce a new function
\begin{equation}
\label{Q_func}
Q(\xi,\eta)=-\frac{1}{2}\int_{\xi/2}^{\eta/2}q(\alpha)\,d\alpha
\end{equation}
and the operator $K:C(\mathbb{R}^2)\mapsto C(\mathbb{R}^2)$ by the
rule
\begin{equation}
\label{K_oper} \Bigl(Kv\Bigr)(\xi,\eta)=\frac{1}{4}
\int_0^\xi\,d\xi_1\int_\xi^\eta\,d\eta_1q\Bigl(\frac{\eta_1-\xi_1}{2}\Bigr)v(\xi_1,\eta_1).
\end{equation}
Rewriting $(\ref{wave_gursa_new_2})$ as
\begin{equation}
\label{Oper_eqn} v=Q-Kv
\end{equation}
and formally solving it by iterations, we get
\begin{equation}
\label{gursa_sol_sum} v(\xi,\eta)=Q(\xi,\eta)+\sum_{n=1}^\infty
(-1)^n(K^nQ)(\xi,\eta).
\end{equation}
To prove the convergence of $(\ref{gursa_sol_sum})$ we need
suitable estimates for $|K^nQ|(\xi,\eta)$. Observe that
\begin{equation}
|Q(\xi,\eta)|\leqslant
\frac{1}{2}\int_0^{\eta/2}|q(\alpha)|\,d\alpha=: S(\eta).
\end{equation}
For the first iteration we have
\begin{eqnarray*}
|(KQ)(\xi,\eta)|\leqslant \frac{1}{4}
\int_0^\xi\,d\xi_1\int_\xi^\eta\,d\eta_1\Bigl|q\Bigl(\frac{\eta_1-\xi_1}{2}\Bigr)\Bigr|S\Bigl(\eta_1\Bigr)\\
\leqslant\frac{S(\eta)}{2}\int_0^\xi\,d\xi_1
\int_{\frac{\xi-\xi_1}{2}}^{\frac{\eta-\xi_1}{2}}|q(\tau)|\,d\tau\leqslant
\frac{S^2(\eta)}{2}\xi.
\end{eqnarray*}
Easy induction argument yields the following estimate
\begin{equation}
\label{gursa_ind} |(K^nQ)(\xi,\eta)|\leqslant
\frac{S^{n+1}(\eta)}{2^n}\frac{\xi^n}{n!},\quad n\in \mathbb{N}.
\end{equation}
Combining (\ref{gursa_sol_sum}) and (\ref{gursa_ind}) one has
\begin{equation}
|v(\xi,\eta)|\leqslant
S(\eta)\operatorname{exp}\Bigl\{S(\eta)\frac{\xi}{2}\Bigr\},
\end{equation}
which due to (\ref{Var_change}) implies $(\ref{gursa_est})$.
Differentiating $(\ref{wave_gursa_new_2})$ we can obtain formulas
for the derivatives of $v$:
\begin{eqnarray}
v_\eta(\xi,\eta)=-\frac{1}{4}q\Bigl(\frac{\eta}{2}\Bigr)-\frac{1}{4}\int_0^\xi
q\Bigl(\frac{\eta-\zeta}{2}\Bigr)
v(\zeta,\eta)\,d\zeta,\label{Deriv_Gur_1}\\
v_\xi(\xi,\eta)=\frac{1}{4}q\Bigl(\frac{\xi}{2}\Bigr)-\frac{1}{4}\int_\xi^\eta
q\Bigl(\frac{\zeta-\xi}{2}\Bigr)v(\xi,\zeta)\,d\zeta+\label{Deriv_Gur_2}\\
\frac{1}{4}\int_0^\xi
q\Bigl(\frac{\xi-\zeta}{2}\Bigr)v(\zeta,\xi)\,d\zeta,\notag
\end{eqnarray}
which lead to $(\ref{gursa_deriv})$.

When $q\in C_{loc}(\mathbb{R}_+)$, we can differentiate one more
time in $(\ref{Deriv_Gur_1})$ with respect to $\xi$ which proves
that $v(\xi,\eta)$ is a classical solution of
$(\ref{gursa_trans})$. For $w(x,t)$ this implies that
$w_{tt}-w_{xx}\in C_{loc}(\mathbb{R}_+\times \mathbb{R}_+)$ and
equation in $(\ref{gursa})$ holds in the classical sense.

When $q\in C^1_{loc}(\mathbb{R}_+)$, we can use
$(\ref{Deriv_Gur_1})$, $(\ref{Deriv_Gur_2})$ to show that $w$ has
all continuous derivatives of the first and second order and thus
is classical.
\end{proof}
Let us emphasize the following simple observation.
\begin{remark}
\label{gursa_rem_uniq} Solution to the boundary value problem
(\ref{gursa}) is unique.
\end{remark}

Let $\{q_n\}_{n=1}^\infty\subset C^\infty(\mathbb{R}_+)$ be such
that
\begin{equation}
q_n\stackrel{L^1_{loc}}\longrightarrow q\quad \text{as
$n\to\infty$};
\end{equation}
by $w_n(x,s)$ we denote the solution of $(\ref{gursa})$
corresponding to the potential $q_n$.

\begin{proposition}
\label{gursa_conv_prop} For solutions $w_n$, $w$ the following
relations hold:
\begin{eqnarray}
w_n\stackrel{C_{loc}}\longrightarrow w,\quad \text{as
$n\to\infty$},\label{gursa_conv}\\
\frac{\partial}{\partial
t}{w_n}\stackrel{L^1_{loc}}\longrightarrow
\frac{\partial}{\partial t}w,
\quad \text{as $n\to\infty$},\label{gursa_deriv_conv_1}\\
\frac{\partial}{\partial x}{w_n}
\stackrel{L^1_{loc}}\longrightarrow \frac{\partial}{\partial x}w,
\quad \text{as $n\to\infty$}\label{gursa_deriv_conv_2}.
\end{eqnarray}
\end{proposition}
\begin{proof}
It is sufficient to prove only $(\ref{gursa_conv})$, since
$(\ref{gursa_deriv_conv_1})$ and $(\ref{gursa_deriv_conv_2})$
follows from that and formulas for derivatives
$(\ref{Deriv_Gur_1})$, $(\ref{Deriv_Gur_2})$. We prove the
convergence for the sequence $\{v_n\}_{n=1}^\infty$, that is
obtained from $\{w_n\}_{n=1}^\infty$ by the change of variables
$(\ref{Var_change})$.

Let us set $\Omega_N=[0,N]\times[0,N]$ and take arbitrary
subsequence from $\{v_n\}_{n=1}^\infty$; we keep the same
notations for it. It is straightforward to check that sequence
$\{v_n\}_{n=1}^\infty$, being restricted to the compact
$\Omega_N$, satisfies the conditions of the Arzela-Ascoli theorem
in $C(\Omega_N)$. Then there exist such a function $\widetilde v\in
C(\Omega_N)$ that for some subsequence
\begin{equation}
v_{n_k} \to \widetilde v,\quad \text{in $C(\Omega_N)$},
\end{equation}
as $k\to\infty$. We rewrite $(\ref{Oper_eqn})$ as
\begin{equation}
v=\mathbf{Q}q-\mathbf{K}(q)v,
\end{equation}
where $\mathbf{Q}: L^1(0,2N)\mapsto C(\Omega_N)$ is defined by
$(\ref{Q_func})$ and $\mathbf{K}:L^1(0,N)\times C(\Omega_N)\mapsto
C(\Omega_N)$ is defined by $(\ref{K_oper})$. We have
\begin{eqnarray}
\label{conv_1}
v-v_{n_k}=\mathbf{Q}(q-q_{n_k})-\mathbf{K}(q-q_{n_k})(v+v_{n_k})-\mathbf{K}(q_{n_k})v+\mathbf{K}(q)v_{n_k}.
\end{eqnarray}
Going to the limit in $(\ref{conv_1})$ we see that
\begin{eqnarray}
\label{conv_2} v-\widetilde
v=-\mathbf{K}(q)v+\mathbf{K}(q)\widetilde v.
\end{eqnarray}
Thus $\widetilde v$ satisfies the equation
\begin{eqnarray}
\label{conv_3} \widetilde v=\mathbf{Q}q-\mathbf{K}(q)\widetilde v.
\end{eqnarray}
Since the solution to $(\ref{Oper_eqn})$ is unique (see remark
\ref{gursa_rem_uniq}), $\widetilde v\equiv v$. Thus, every
subsequence of $\{v_n\}_{n=1}^\infty$ contains a subsequence,
convergent to $v$ in $C(\Omega_N)$. It implies that the very
sequence converges to $v$. Since $N$ is arbitrary, we arrive at
$(\ref{gursa_conv})$.
\end{proof}

\section{Acknowledgments}

The results of this paper were discussed with M.~Belishev,
S.~Ivanov, K. Mirzoev, A.~Rybkin, and A.~Shkalikov.  The authors
are very grateful to all of them. These results were presented at
the International Meeting on Inverse and Spectral Problems at the
University of Auckland (December, 2007) and at the Seminar on
Mathematical Physics at the California Institute of Technology
(February, 2008). The authors are  grateful to participants of the
both meetings for fruitful discussions.

Research of Sergei Avdonin was supported in part by the NSF, grant
ARC 0724860.

\end{document}